\documentclass[11pt]{article}
\usepackage{multicol}
\usepackage{epic}
\usepackage{epsfig}
\usepackage[english]{babel}

\usepackage{amssymb,epsf,amsmath}
\newtheorem{Theorem}{Theorem}[section]
\newtheorem{Proposition}[Theorem]{Proposition}

\newtheorem{Lemma}[Theorem]{Lemma}

\newtheorem{Corollary}[Theorem]{Corollary}

\newtheorem{Ex}{Example}[section]
\newenvironment{Example}{\begin{Ex}\rm}{\smallskip\end{Ex}}

\newtheorem{proofhead}{Proof}

\newcommand{\qed}{\hfill$\Box$}

\newcounter{mylistcnt}
\renewcommand{\themylistcnt}{{\rm({\roman{mylistcnt}})}}

\newcounter{zad}

\begin{document}
\title{Pickands-Piterbarg constants for self-similar Gaussian processes}
\author{
Krzysztof D\c{e}bicki$^{1}$, Kamil Tabi\'{s}$^{1}$\\
%\thanks{KD work was supported by NCN Research
%Grant No 2011/01/B/ST1/01521 (2011-2013)} ,
%\thanks{KT work was supported by NCN Grant ...}
\emph{ $^1$ Mathematical Institute, University of Wroc\l aw, pl. Grunwaldzki 2/4, 50-384 Wroc\l aw, Poland}
}
\date{}
\maketitle

\begin{abstract}
For a centered self-similar
Gaussian process $\{Y(t):t\in[0,\infty)\}$ and $R\ge0$ we analyze asymptotic behaviour of
\[
\mathcal{H}_Y^R(T) \; = \; \mathbf{E} \exp  \left( \sup_{t \in [0,T]} \sqrt{2} Y(t) - (1+R) \sigma_Y^2(t) \right),
\]
as $T\to\infty$.
We prove that
$\mathcal{H}_Y^R=\lim_{T\to\infty} \mathcal{H}_Y^R(T)\in(0,\infty)$, for $R>0$
and
\[\mathcal{H}_Y=\lim_{T\to\infty} \frac{\mathcal{H}_Y^0(T)}{T^\gamma}\in(0,\infty)\]
for suitably chosen $\gamma>0$.
Additionally, we find bounds for $\mathcal{H}_Y^R$, $R>0$ and
a surprising relation between $\mathcal{H}_Y$
and
classical Pickands constants.
\\


\end{abstract}

%\newpage
%%%%%%%%%%%%%%%%%%%%%%%%%%%%%%%%%%%%%%%%%%%%%%%%%%%%%%%%%%%%%%%%%%%%%%%%%%%%%%%%%%%%%%%%%%%%%%%%%%%%%%%%

\section{Introduction}\label{s.intro}

For a centered Gaussian process
$ \{ Y(t) : t \in [0,\infty) \} $
with a.s. continuous sample paths,
$\mathbf{Var}(Y(t))=\sigma^2_Y(t)$ and $ Y(0) = 0 $ a.s.,
let
\begin{eqnarray}
\mathcal{H}_Y^R(T) \; = \; \mathbf{E} \exp  \left( \sup_{t \in [0,T]} \sqrt{2} Y(t) - (1+R) \sigma_Y^2(t) \right),
\label{wills}
\end{eqnarray}
where $R\ge0$ and let
$\mathcal{H}_Y(T):=\mathcal{H}_Y^0(T)$.

Functionals $\mathcal{H}_Y^R(T),\mathcal{H}_Y(T)$
play important role in many areas of probability theory.
For example, consider a fractional Brownian
motion  $\{B_\kappa(t):t\in[0,\infty)\}$ with Hurst parameter $\kappa/2\in(0,1]$, i.e.
a centered Gaussian process with stationary increments, continuous sample paths a.s and
variance function $\mathbf{Var}(B_\kappa(t))=t^\kappa$.
Then, for $\kappa\in(0,2]$, {\it Pickands constants} $ \mathcal{H}_{B_\kappa} $ defined as
\begin{eqnarray}
\mathcal{H}_{B_\kappa} \; = \; \lim_{T \to \infty} \frac{\mathcal{H}_{B_\kappa}(T)}{T} , \label{def.p1}
\end{eqnarray}
and  {\it Piterbarg constants}
$ \mathcal{H}_{B_\kappa}^R $, for $ R > 0 $, defined as
\begin{eqnarray}
\mathcal{H}_{B_\kappa}^R \; = \; \lim_{T \to \infty} \mathcal{H}_{B_\kappa}^R(T)  \label{def.p2}
\end{eqnarray}
 play a key role in the
extreme value theory of Gaussian processes; see, e.g., \cite{Pic69, Pit96, PiP79} or more
recent contributions \cite{HKP19,Mic17}.
In \cite{Deb02} it was observed that the notion of Pickands and Piterbarg constants
can be extended to  {\it generalized} Pickands and Piterbarg constants, defined as
\[\mathcal{H}_{\eta} \; = \; \lim_{T \to \infty} \frac{\mathcal{H}_{\eta}(T)}{T}\]
and
\[\mathcal{H}_{\eta}^R \; = \; \lim_{T \to \infty} \mathcal{H}_{\eta}^R(T)\]
respectively, where $R>0$ and $\{\eta(t):t\in(0,\infty)\}$ is a centered Gaussian process with stationary increments.
We refer to \cite{Albin,DHJ18,Mic04, DHL17,DMR03,DiM15}
for properties and other representations of $\mathcal{H}_{B_\kappa}$, $\mathcal{H}_{B_\kappa}^R$
and generalized Pickands-Piterbarg constants,
and to \cite{DHJR18,DHT1} for multidimensional analogs of Pickands-Piterbarg constants.

Recently, see e.g. \cite{DHL17}, it was found that
for general Gaussian processes $Y$ (satisfying some regularity conditions)
functionals (\ref{wills}) appear in the formulas for exact asymptotics of supremum of some
Gaussian processes; see Proposition \ref{Cor09}.

The interest in analysis of properties of (\ref{wills}) stems also from
an important contribution
\cite{DiM15} which established a direct connection between
Pickands constants and max-stationary stable processes; see also \cite{DEH17, DeH16, DeH19}.

Constants $\mathcal{H}_Y(T)$ appear also in the context of convex geometry
where they are known  as {\it Wills functionals};
see \cite{Wit96}.

In this contribution we analyze properties of
$\mathcal{H}_Y^R(T)$ and $\mathcal{H}_Y(T)$
for the class of general self-similar Gaussian processes $Y$ with non stationary increments.
In particular, we find  analogs of limits (\ref{def.p1}), (\ref{def.p2})
and give some bounds for them.
Surprisingly, it appears that, up to some explicitly given constant,
$\mathcal{H}_Y$ is equal
to the classical  $\mathcal{H}_{B_\kappa}$ for some appropriately chosen $\kappa$.

\section{Notation and preliminary results}\label{s.notation}
%\section{An extension of Pickands Lemma} \label{S32}
Let $ \{ Y(t) : t \geq 0 \} $ be a centered Gaussian process with a.s. continuous sample paths and let
\[V_{Y}(s,t) := \mathbf{Var}(Y(s) - Y(t)) ; \]
\[ R_{Y}(s,t):= \mathbf{Cov}(Y(s),Y(t)) \; . \]

We say that a stochastic process $ Y(\cdot) $ is {\it self-similar} with index $ H > 0 $, if for all $ a > 0 $,
\begin{equation}
\label{lab97} \left\{ Y(at) : t \geq 0 \right\} \; \stackrel{\mathcal{D}}{=} \; \left\{ a^{H} Y(t) : t \geq 0 \right\} \; .
\end{equation}
A straightforward consequence of (\ref{lab97}) is that for self-similar Gaussian processes,
$ \sigma_{Y}^2(t) = \sigma_{Y}^2(1) t^{2H} $ for $ t \geq 0 $.

We write $ Y \in \mathbf{S}(\alpha,\kappa,c_Y) $ if
\begin{itemize}
\item[\textbf{S1}] $ Y(\cdot) $ is self-similar with index $ \alpha/2 > 0 $ and $ \sigma_{Y}^2(1) = 1 $;
\item[\textbf{S2}] there exist $ \kappa \in (0,2] $ and $ c_{Y} > 0 $ such that $ \mathbf{Var}(Y(1)-Y(1-h)) = c_Y |h|^\kappa + o(|h|^\kappa) $, as $ h \to 0 $. \\
\end{itemize}

It is well known (see Lamperti \cite{Lam}) that $ \{ Y(t) : t \geq 0 \} $
is a self-similar Gaussian process with index
$ \alpha/2 $, if and only if, its Lamperti transform $ X(t) = e^{-(\alpha / 2)t} Y(e^{t}) $
is a stationary Gaussian process.
Thus, there is a unique correspondence between self-similar Gaussian processes and stationary
Gaussian processes. In fact condition \textbf{S2} relates to regularity condition for a
covariance function of the stationary counterpart of $Y$.
More precisely, let $ \{ X(t) : t \in \mathbb{R} \} $
be a stationary Gaussian process such that $ R_X(t,0) = 1 - a|t|^{\kappa} + o(|t|^{\kappa}) $, as $ t \to 0 $,
for $ \kappa \in (0,2] , a > 0 $. Then, one can check that
the self-similar process  $ Y(t) := t^{\alpha/2} X(\log t) $, for $ \alpha \in (0,2] $,
is
$\mathbf{S}(\alpha,\kappa,c_Y) $ with
\[ c_{Y} \; = \; \left\{
\begin{array}{lcr} 2 a & \textrm{for} & \kappa < 2 \\[1ex] \frac{\alpha^2}{4} + 2 a & \textrm{for} & \kappa = 2 \end{array} \right. .
\]

Below we specify some important classes of self-similar Gaussian processes that
satisfy
\textbf{S1-S2}.\\
$\diamond$
\emph{Fractional Brownian motion} $ B_{\alpha}\in \mathbf{S}(\alpha,\alpha,1)$ with $ \alpha/2 \in (0,1] $.\\
$\diamond$
\emph{Bifractional Brownian motion} $ \{ Y^{(1)}(t) : t \geq 0 \} $ with parameters $ \alpha \in (0,2) $ and $ K \in (0,1] $ is a centered Gaussian process with covariance function
\[ R_{Y^{(1)}}(t,s) \; = \; \frac{1}{2^K}  \left( (t^\alpha + s^\alpha)^K - |t-s|^{\alpha K} \right),  \]
see e.g. \cite{Ho},\cite{Lei}.
We have $ Y^{(1)} \in \mathbf{S}(\alpha K,\alpha K,2^{1-K}) $.\\
$\diamond$
\emph{Sub-fractional Brownian motion}  $ \{ Y^{(2)}(t) : t \geq 0 \}  $ with parameter $ \alpha \in (0,2) $ is a centered Gaussian process
with covariance function
\[ R_{Y^{(2)}}(t,s) \; = \; \frac{1}{2-2^{\alpha-1}}  \left( t^\alpha + s^\alpha - \frac{(t+s)^\alpha + |t-s|^{\alpha}}{2} \right), \]
see \cite{Bojdecki,DzZ04}.
Then $ Y^{(2)} \in \mathbf{S}(\alpha,\alpha,(2-2^{\alpha-1})^{-1}) $.
\\
$\diamond$
\emph{$k$-fold integrated fractional Brownian motion} $\{ Y^{(3),k}(t) : t \geq 0 \}  $
with parameters $ k \in \mathbb{N}:=\{1,2,3,...\} $ and $ \alpha \in (0,2] $ is a Gaussian process defined as
\begin{eqnarray*}
Y^{(3),1}(t) \; = \; \sqrt{\alpha+2} \int_0^t B_\alpha (s) ds & \quad \textrm{for} & k=1 \; ; \\[1ex]
Y^{(3),k}(t) \; = \; \sqrt{\frac{k(\alpha+2k)(\alpha+k-1)}{\alpha+2k-2}} \int_0^t Y^{(3),k-1}(s) ds & \quad \textrm{for} & k \geq 2 \; . \\
\end{eqnarray*}
Then
 $Y^{(3),k} \in \mathbf{S}(\alpha+2k,2,\frac{k(\alpha+2k)(\alpha+k-1)}{\alpha+2k-2}) $ for $ k \geq 1 $. \\
$\diamond$
\emph{Time-average of fractional Brownian motion} $\{ Y^{(4)}(t) : t \geq 0 \}  $
with parameter $ \alpha \in (0,2] $ is a Gaussian process defined as
\[ Y^{(4)}(t) \; = \; \sqrt{\alpha+2} \: \frac{1}{t} \int_0^t B_\alpha (s)ds \; . \]
Its covariance function is of the form
\[ R_{Y^{(4)}}(t,s) \; = \; \frac{ (\alpha+2)  \left(s^{\alpha+1}t + st^{\alpha+1}\right) + |t-s|^{\alpha+2} - t^{\alpha+2} - s^{\alpha+2} }{2(\alpha+1)ts} \;  \]
and we have $ Y^{(4)} \in \mathbf{S}(\alpha,2,1) $.
\\
$\diamond$
\emph{Dual fractional Brownian motion} $\{ Y^{(5)}(t) : t \geq 0 \}$ with parameter  $ \alpha \in (0,2] $ is a centered Gaussian process
defined as
\[ Y^{(5)}(t) \; = \; t^{\alpha+1} \sqrt{\frac{2}{\Gamma(\alpha+1)}} \int_0^{\infty} B_\alpha(s) e^{-st}ds \; , \]
see \cite{LiShao}.
We have
\[ R_{Y^{(5)}}(t,s) \; = \; \frac{t^{\alpha}s + s^{\alpha}t}{t+s} \;  \]
and $ Y^{(5)} \in \mathbf{S}(\alpha,2,\alpha/2) $.
\\
\\

In the rest of the  paper, let $\overline{X}(s):=X(s)/\sigma_X(s)$.
$\Psi(\cdot)$ denotes the tail distribution function of the standard normal random variable.

The following proposition plays key role in the proofs of the main results of this contribution,
confirming also that
functionals
$\mathcal{H}_{Y}(\cdot)$ and $\mathcal{H}_{Y}^R(\cdot)$
for  $ Y \in \mathbf{S}(\alpha,\kappa,c_Y) $  appear in
the asymptotics of extremes of Gaussian processes.

\begin{Proposition} \label{Cor09}
Let $ Y \in \mathbf{S}(\alpha,\kappa,c_Y) $ and let $ \{ X(t) : t \geq 0 \} $ be a centered Gaussian process  with
$ R_{\overline{X}}(t,s) = \exp(-aV_Y(t,s)) $ for $ a > 0 $ and $ \sigma_X(t) = \frac{1}{1+bt^{\beta}} $
for $ b \geq 0, \beta>0 $.
%, $ \beta \geq \lambda $, where $ \lambda \in (0,2] $ is given. Let
%\[ R \; = \; \left\{ \begin{array}{rcl} b/a & \textrm{if} & \beta = \lambda \\[1ex] 0 & \textrm{if} & \beta > \lambda \end{array} \right. . \]
\begin{itemize}
\item[\emph{(i)}] If $ \alpha = \beta $, then as $ u \to \infty $,
\[ \mathbf{P}  \left( \sup_{t \in [0,Tu^{-2/\alpha}]} X(t) > u \right) \; = \; \mathcal{H}_{Y}^{b/a}\! \left(a^{1/\alpha}T\right) \Psi(u) (1+o(1)) \; . \]
\item[\emph{(ii)}] If $ \alpha < \beta $, then as $ u \to \infty $,
\[ \mathbf{P}  \left( \sup_{t \in [0,Tu^{-2/\alpha}]} X(t) > u \right) \; = \; \mathcal{H}_{Y}\! \left(a^{1/\alpha}T\right)\Psi(u) (1+o(1)) \; . \] \vspace{2ex}
\end{itemize}
\end{Proposition}
The proof of Proposition \ref{Cor09} is given in Section \ref{s.PrPr1}.

\section{Pickands-Piterbarg constants for self-similar Gaussian processes} \label{S33}

The aim of this section is to find analogs of
Pickands and Piterbarg constants for the class of
self-similar Gaussian processes $ Y \in \mathbf{S}(\alpha,\kappa,c_Y) $.

\subsection{Piterbarg constants}
For $R>0$ and $ Y \in \mathbf{S}(\alpha,\kappa,c_Y) $ let us introduce an
analog of Piterbarg constant $\mathcal{H}_{B_\kappa}^R$ as follows
\[
\mathcal{H}_{Y}^R  :=  \lim_{T \to \infty} \mathcal{H}_{Y}^R(T) \; = \;
\lim_{T \to \infty} \mathbf{E} \exp  \left( \sup_{t \in [0,T]} \sqrt{2} Y(t) - (1+R) t^\alpha \right) .
\]
In the next theorem we prove  that $\mathcal{H}_{Y}^R $ is well-defined and
compare $\mathcal{H}_{Y}^R$ with classical Piterbarg constants.

\begin{Theorem} \label{TwDr2}
Let $ Y \in \mathbf{S}(\alpha,\kappa,c_Y) $. Then, for any $ R > 0 $,
\[ \mathcal{H}_Y^R \in(0,\infty).\]
Furthermore
\[ \mathcal{H}_{B_\kappa}^{R/c_1} \; \leq \; \mathcal{H}_Y^R \; \leq \; \mathcal{H}_{B_\kappa}^{R/c_2} \; , \]
where
\[ c_1 = \inf_{x \in [0,1)} \frac{V_Y(1,x^{\kappa/\alpha})}{|1-x|^\kappa} \qquad \textrm{and} \qquad c_2 = \sup_{x \in [0,1)} \frac{V_Y(1,x^{\kappa/\alpha})}{|1-x|^\kappa} \;. \]
\end{Theorem}
The proof of Theorem \ref{TwDr2}  is  given in Section \ref{PrTwDr2}.
%The following lower bound holds.
\begin{Proposition} \label{Cor11}
Let $ Y \in \mathbf{S}(\alpha,\kappa,c_Y) $. Then
\[ \mathcal{H}_Y^R \; \geq \; \frac{1}{2}  \left( 1 + \sqrt{1 + \frac{1}{R}} \right) . \]
\end{Proposition}
The proof of Proposition \ref{Cor11} is postponed to Section \ref{PrCor11}.

The following corollary follows from Theorem  \ref{TwDr1}
combined with the fact that
$\mathcal{H}_{B_1}^R \; = \; 1 + \frac{1}{R}$ (see, e.g., \cite{DebMan})
and
$\mathcal{H}_{B_2}^R \; = \; \frac{1}{2}  \left( 1 + \sqrt{1 + \frac{1}{R}} \right)$
(see, e.g., \cite{KonsPit})

\begin{Corollary} \label{Cor13}  Let $ Y \in \mathbf{S}(\alpha,\kappa,c_Y) $.
\\
(i) If $\kappa=1$, then
\[ 1 + \frac{1}{R}  \left(\inf_{x \in [0,1)} \frac{V_Y(1,x^{1/\alpha})}{|1-x|}\right) \; \leq \; \mathcal{H}_{Y}^R \; \leq \; 1 + \frac{1}{R}  \left(\sup_{x \in [0,1)} \frac{V_Y(1,x^{1/\alpha})}{|1-x|}\right) . \] \vspace{2ex}
(ii) If $\kappa=2$, then
\[ \frac{1}{2}  \left( 1 + \sqrt{1+\frac{1}{R}} \right) \; \leq \; \mathcal{H}_{Y}^R \; \leq \; \frac{1}{2}
\left( 1 + \sqrt{1+\frac{1}{R} \left(\sup_{x \in [0,1)} \frac{V_Y(1,x^{2/\alpha})}{|1-x|^2}\right)} \right) . \]
\end{Corollary}

In the following example we specify Corollary \ref{Cor13}
for some particular self-similar processes introduced in Section \ref{s.notation}.

\begin{Example}
The following bounds hold.\\
$\diamond$
$k-$fold integrated fractional Brownian motion $Y^{(3),k}$:
\[ \frac{1}{2}  \left( 1 + \sqrt{1+\frac{1}{R}} \right) \; \leq \;
\mathcal{H}_{Y^{(3),k}}^R \; \leq \; \frac{1}{2}
\left( 1 + \sqrt{1+ \frac{4k(\alpha+k-1)}{R(\alpha+2k)(\alpha+2k-2)}} \right).\]
$\diamond$
Time-average of fractional Brownian motion $Y^{(4)}$ with parameter $ \alpha \in (0,2] $:
\[ \frac{1}{2}  \left( 1 + \sqrt{1+\frac{1}{R}} \right) \; \leq \; \mathcal{H}_{Y^{(4)}}^R \;
\leq \; \frac{1}{2}  \left( 1 + \sqrt{1+ \frac{4}{R(\alpha+2)}} \right) . \]
The above bounds improve results obtained in \cite{DebTab}  for constants
\begin{eqnarray*}
\mathcal{F}_{\alpha} & = & \lim_{T \to \infty} \mathbf{E} \exp  \left( \sup_{t \in (0,T]} \frac{1}{t} \int_0^t  \left( \sqrt{2} B_{\alpha}(s) - s^{\alpha} \right)  ds \right) \\[1ex]
\nonumber & = & \lim_{T \to \infty} \mathbf{E} \exp  \left( \sup_{t \in (0,T]} Y^{(4)}(t) - \frac{\alpha+2}{\alpha+1} t^{\alpha}
\right) \; = \; \mathcal{H}_{Y^{(4)}}^{1/(\alpha+1)} \; ,
\end{eqnarray*}
leading to
\[ \frac{1}{2}  \left( 1 + \sqrt{2 + \alpha} \right) \; \leq\; \mathcal{F}_{\alpha} \; \leq \; \frac{1}{2}  \left( 1 + \sqrt{1 + 4(\alpha+1)/\alpha^2} \right),  \]
while in \cite{DebTab} it was proved that
$ \mathcal{F}_{\alpha} \leq 2 + \alpha $ for $ \alpha \in [1,2) $.\\
$\diamond$
Dual fractional Brownian motion $Y^{(5)}$ with parameter $\alpha\in(0,2]$:
\[ \frac{1}{2}  \left( 1 + \sqrt{1+\frac{1}{R}} \right) \; \leq \; \mathcal{H}_{Y^{(5)}}^R \; \leq \; \frac{1}{2}  \left( 1 + \sqrt{1+\frac{2}{R\alpha}} \right).  \]
\end{Example}

\subsection{Pickands constants}
In this section we focus on an analog of Pickands constants
for $ Y \in \mathbf{S}(\alpha,\kappa,c_Y) $.
Let
\begin{eqnarray*}
\mathcal{H}_{Y} & := & \lim_{T \to \infty} \frac{\mathcal{H}_{Y}(T)}{T^{\alpha/\kappa}} \; = \; \lim_{T \to \infty} \frac{\mathbf{E} \exp  \left( \sup_{t \in [0,T]} \sqrt{2} Y(t) - t^\alpha \right)}{T^{\alpha/\kappa}} \;.
\end{eqnarray*}
We observe that for $Y(t)=B_{\kappa}(t)$ the above definition agrees with
the notion of
the classical Pickands constant $\mathcal{H}_{B_{\kappa}}$, since $\alpha=\kappa$ in this case.

In the following theorem we show that $\mathcal{H}_{Y}$ is well-defined
and find a surprising relation between $\mathcal{H}_{Y}$ and $\mathcal{H}_{B_{\kappa}}$.

\begin{Theorem} \label{TwDr1}
Let $ Y \in \mathbf{S}(\alpha,\kappa,c_Y) $. Then $ \mathcal{H}_{Y}\in(0,\infty)$ and
\[ \mathcal{H}_{Y} \; = \; \frac{\kappa}{\alpha}  \left( c_Y \right)^{1/\kappa} \mathcal{H}_{B_\kappa} \; . \]
\end{Theorem}

%\begin{Remark} \label{Rem21} \emph{
%Note that for $ Y \in \mathbf{S}(\alpha,\kappa,c_Y) $, in general there is no property of subadditivity of $ \mathcal{H}_{Y}(\cdot) $ similar to one in case Pickands or generalized Pickands constants. It suffices to consider $ Y \in \mathbf{S}(\alpha,\kappa,c_Y) $ with $ \alpha > \kappa $. \\
%} \end{Remark}

Complete proof of Theorem \ref{TwDr1} is presented in Section \ref{PrTwDr1}.
\\

The following corollary is an immediate consequence of Theorem  \ref{TwDr1}
and the fact that $ \mathcal{H}_{B_1}=1 $ and $ \mathcal{H}_{B_2}=\frac{1}{\sqrt{\pi}} $.

\begin{Corollary} \label{Cor07} Let $ Y \in \mathbf{S}(\alpha,\kappa,c_Y) $.
\\
(i) If $\kappa=1$, then
\[ \mathcal{H}_{Y} \; = \; \lim_{T \to \infty} \frac{\mathbf{E} \exp  \left( \sup_{t \in [0,T]} \sqrt{2} Y(t) - t^\alpha \right)}{T^{\alpha}} \; = \; \frac{c_{Y}}{\alpha} \; .\]
(ii) If $\kappa=2$, then
\[ \mathcal{H}_{Y} \; = \; \lim_{T \to \infty} \frac{\mathbf{E} \exp  \left( \sup_{t \in [0,T]} \sqrt{2} Y(t) - t^\alpha \right)}{T^{\alpha/2}} \; = \; \frac{2}{\alpha} \sqrt{\frac{c_{Y}}{\pi}} \;.\]
\end{Corollary}

In the following example we specify the findings of this section to
self-similar Gaussian processes introduced in Section \ref{s.notation}.
\begin{Example}
The following equalities hold.\\
$\diamond$ Bifractional Brownian motion with parameters $ \alpha \in (0,2) $ and $ K \in (0,1] $:
\[ \mathcal{H}_{Y^{(1)}} \; = \; \lim_{T \to \infty} \frac{\mathbf{E} \exp  \left( \sup_{t \in [0,T]} \sqrt{2} Y^{(1)}(t) - t^{\alpha K} \right)}{T} \; = \; 2^{\frac{1-K}{\alpha K}} \mathcal{H}_{B_{\alpha K}} \; . \]
$\diamond$
Sub-fractional Brownian motion with parameter $ \alpha \in (0,2) $:
\[ \mathcal{H}_{Y^{(2)}} \; = \; \lim_{T \to \infty} \frac{\mathbf{E} \exp  \left( \sup_{t \in [0,T]} \sqrt{2} Y^{(2)}(t) - t^\alpha \right)}{T} \; = \; (2-2^{\alpha-1})^{-1/\alpha} \mathcal{H}_{B_{\alpha}} \; . \]
$\diamond$
k-fold integrated fractional Brownian motion with parameters $ k \in \mathbb{N} $ and $ \alpha \in (0,2] $:
\[ \mathcal{H}_{Y^{(3),k}} =  \lim_{T \to \infty} \frac{\mathbf{E} \exp  \left( \sup_{t \in [0,T]} \sqrt{2} Y^{(3),k}(t) - t^{k+\alpha/2} \right)}{T^{k+\alpha/2}} \; = \; \sqrt{\frac{4k(\alpha+k-1)}{\pi(\alpha+2k)(\alpha+2k-2)}}  . \]
$\diamond$
Time-average of fractional Brownian motion with parameter $ \alpha \in (0,2] $:
\[ \mathcal{H}_{Y^{(4)}} =  \lim_{T \to \infty}
\frac{\mathbf{E} \exp  \left( \sup_{t \in [0,T]} \sqrt{2} Y^{(4)}(t) - t^\alpha \right)}{T^{\alpha/2}} =
\frac{2}{\sqrt{\pi} \alpha} . \]
$\diamond$
Dual fractional Brownian motion with parameter $ \alpha \in (0,2] $:
\[ \mathcal{H}_{Y^{(5)}} \; = \; \lim_{T \to \infty}
\frac{\mathbf{E} \exp  \left( \sup_{t \in [0,T]} \sqrt{2} Y^{(5)}(t) - t^\alpha \right)}{T^{\alpha/2}} \; = \; \sqrt{\frac{2}{\pi \alpha}} \;. \]
 \end{Example}

\section{Proofs}
In the rest of the paper we use the following notation $v_Y(t):=V_Y(1,t)$.
We begin with the following lemma, skipping its straightforward proof.
\begin{Lemma} \label{Lem23}
Let $ Y \in \mathbf{S}(\alpha,\kappa,c_Y) $ and $ \hat{Y}_{\alpha_1} = Y(t^{\alpha_1}) $, for some $ \alpha_1 > 0 $. Then, for any $ R \geq 0 $ and $ T > 0 $,
\begin{itemize}
\item[\emph{(i)}] $ \mathcal{H}_{cY}^R(T) \; = \; \mathcal{H}_{Y}^R \! \left(c^{2/\alpha}T\right) $ for any $ c > 0 $;
\item[\emph{(ii)}] $ \mathcal{H}_{\hat{Y}_{\alpha_1}}^R(T) \; = \; \mathcal{H}_{Y}^R(T^{\alpha_1}) $. \\
\end{itemize}
\end{Lemma}

\subsection{Proof of Proposition \ref{Cor09}} \label{s.PrPr1}
In the next lemma we present a useful bound on $ V_Y(\cdot,\cdot) $, for $ Y \in \mathbf{S}(\alpha,\kappa,c_Y) $.

\begin{Lemma} \label{Lem10}
Let $ Y \in \mathbf{S}(\alpha,\kappa,c_Y) $. Then there exists positive constant $ C $,
such that for $ \gamma = \min(\alpha,\kappa) $,
$ T > 0 $ and all $ t,s \in [0,T] $,
\[  V_{Y}(t,s) \; \leq \; C T^{\alpha-\gamma} |t-s|^{\gamma} \; . \]
\end{Lemma}
\noindent \textbf{Proof:}
For $ t = s $ the thesis is obvious. Suppose that $ 0 \leq s < t \leq T $ and let $ \epsilon \in (0,1) $ be such that
for $ \delta \in (0,1) $
\[ (1-\epsilon) c_{Y} |1-x|^\kappa \; \leq \; V_{Y}(1,x) \; \leq \; (1+\epsilon) c_{Y} |1-x|^\kappa \; , \]
for all $ x \in [\delta,1] $ (due to \textbf{S2}). For $ s/t \geq \delta $ we have
\begin{eqnarray*}
V_{Y}(t,s) & = & t^\alpha V_{Y}(1,s/t) \; \leq \; t^\alpha (1+\epsilon) c_{Y} |1-s/t|^\kappa \; \leq \; t^\alpha (1+\epsilon) c_{Y} |1-s/t|^{\min(\alpha,\kappa)} \\[1ex]
& = & t^{\alpha-\gamma} (1+\epsilon) c_{Y} |t-s|^{\gamma} \; \leq \; T^{\alpha-\gamma} (1+\epsilon) c_{Y} |t-s|^{\gamma}.
\end{eqnarray*}
For $ s/t \leq \delta $ we have  $ |1-\delta|^{\gamma} \leq |1-s/t|^{\gamma} $.
Hence, $ t^{\gamma}|1-\delta|^{\gamma} \leq |t-s|^{\gamma} $.
Then
\[ V_{Y}(t,s) \; = \; t^\alpha V_{Y}(1,s/t) \; \leq \; t^{\gamma} t^{\alpha-\gamma} \frac{|1-\delta|^{\gamma_2}}{|1-\delta|^{\gamma}}
 \max_{x \in [0,\delta]} V_{Y}(1,x) \; \leq \; T^{\alpha-\gamma} \frac{\max_{x \in [0,\delta]} V_{Y}(1,x)}{|1-\delta|^{\gamma}} |t-s|^{\gamma}. \]

Hence the proof is completed with
${C = \max  \left( (1+\epsilon)c_{Y}, \frac{\max_{x \in [0,\delta]} V_{Y}(1,x)}{|1-\delta|^{\gamma}} \right)}$.
\hfill \raisebox{-1.5ex}{$ \square $}
\\
\noindent \textbf{Proof of Proposition \ref{Cor09}:}\\
Using that for any Gaussian process $ Y(\cdot) $, the variogram function $ V_Y(\cdot,\cdot) $ is negative definite,
by Schoenberg theorem, function $ \exp(-V_Y(\cdot,\cdot)) $ is positive definite. Thus there exists a Gaussian process
$ \{ {X}(t) : t \geq 0 \} $ with $ R_{{X}}(t,s) = \exp(-V_Y(t,s)) $.
\\
The rest of the proof follows straightforwardly from  Theorem 2.1 in \cite{DHL17} and Lemma \ref{Lem10}
applied to $X_u(t)=X(tu^{-2/\alpha})$.
\hfill \raisebox{-1.5ex}{$ \square $}

%In this section we proof main results of this contribution.
\subsection{Proof of Theorem \ref{TwDr2}} \label{PrTwDr2}
Let $ R , T > 0 $. Consider $ \hat{Y}(t) = Y(t^{\kappa/\alpha}) $.
Then $ \hat{Y} \in \mathbf{S}(\kappa,\kappa,c_Y (\kappa/\alpha)^\kappa) $ with
$ v_{\hat{Y}}(x) = v_{Y}(x^{\kappa/\alpha}) $ and
\[ c_1 = \inf_{x \in [0,1)} \frac{v_Y(x^{\kappa/\alpha})}{|1-x|^\kappa}\in(0,\infty)\]
\[ c_2 = \sup_{x \in [0,1)} \frac{v_Y(x^{\kappa/\alpha})}{|1-x|^\kappa}\in(0,\infty) \;. \]
 Let $ \{ X(t) : t \geq 0 \} $, $ \{ X_1(t) : t \geq 0 \} $, $ \{ X_2(t) : t \geq 0 \} $ be centered Gaussian processes
 with $ R_{\overline{X}}(t,s) = \exp(-V_{\hat{Y}}(t,s)) $, $ R_{\overline{X}_i}(t,s) = \exp(-c_i V_{B_\kappa}(t,s)) $
 and $ \sigma_X(t) = \sigma_{X_i}(t) = \frac{1}{1+Rt^{\kappa}} $, respectively.
 Then, for all $ t,s \geq 0 $,
\begin{eqnarray*}
R_{\overline{X}_1}(t,s) & = & \exp(-c_1 V_{B_\kappa}(t,s)) \; \geq \; V_{\hat{Y}}(t,s) \\[1ex]
& = & \exp(-V_{\hat{Y}}(t,s)) \; \geq \; \exp(-c_2 V_{B_\kappa}(t,s)) \; = \; R_{\overline{X}_2}(t,s) \;
\end{eqnarray*}
and hence, due to Slepian's inequality (see, e.g., Corollary 2.4 in \cite{Adl90}) we obtain that for all $ u > 0 $,
\[ \mathbf{P}  \left( \sup_{t \in [0,Tu^{-2/\kappa}]} X_1(t) > u \right) \; \leq \; \mathbf{P}  \left( \sup_{t \in [0,Tu^{-2/\kappa}]} X(t) > u \right) \; \leq \; \mathbf{P}  \left( \sup_{t \in [0,Tu^{-2/\kappa}]} X_2(t) > u \right) . \]
Application of Proposition \ref{Cor09} (i) to the inequalities above gives
\begin{equation} \label{lab111} \mathcal{H}_{B_\kappa}^{R/c_1} \! \left(c_1^{1/\kappa}T\right) \; \leq \; \mathcal{H}_{\hat{Y}}^R(T) \; = \; \mathcal{H}_{Y}^R (T^{\kappa/\alpha}) \; \leq \; \mathcal{H}_{B_\kappa}^{R/c_2} \! \left(c_2^{1/\kappa}T\right) , \end{equation}
where equality above follows from Lemma \ref{Lem23} (ii).
Note that all functions in (\ref{lab111}) are increasing and hence,
after sending $ T \to \infty $ in (\ref{lab111}), the proof is completed.
\hfill \raisebox{-1.5ex}{$ \square $}
\subsection{Proof of Proposition \ref{Cor11}} \label{PrCor11}
Since $ \sigma_{Y}^2(t) = t^\alpha $ then the Schwarz inequality implies that $ R_Y(t,s) \leq (ts)^{\alpha / 2} $ for all $ t,s \geq 0 $. Therefore for $ t,s \geq 0 $,
\[ V_{Y}(t,s) \; \geq \; t^{\alpha} + s^{\alpha} - 2(ts)^{\alpha / 2} \; = \; |t^{\alpha / 2}-s^{\alpha / 2}|^2 \; = V_{\hat{B}_2}(t,s) \; , \]
where $ \hat{B}_{2}(t) = B_2(t^{\alpha/2}) $. Thus, due to Slepian's inequality,  we conclude that
\begin{eqnarray*}
\mathcal{H}_{Y}^R(T) & = & \int_{\mathbb{R}} e^x \mathbf{P}  \left( \sup_{t \in [0,T]} \sqrt{2} Y(t) - (1+R)t^{\alpha} > x \right)  dx \\[1ex]
& \geq & \int_{\mathbb{R}} e^x \mathbf{P}  \left( \sup_{t \in [0,T]} \sqrt{2} \hat{B}_{2}(t) - (1+R)t^{\alpha} > x \right)  dx \; = \\[1ex]
& = & \mathcal{H}_{\hat{B}_{2}}^R(T) \; = \; \mathcal{H}_{B_{2}}^R \! \left(T^{\alpha/2}\right) \; ,
\end{eqnarray*}
where the last equality follows from Lemma \ref{Lem23} (ii).
We obtain the thesis after sending $ T \to \infty $ in the inequality above. %This completes the proof.
\hfill \raisebox{-1.5ex}{$ \square $}

\subsection{Proof of Theorem \ref{TwDr1}} \label{PrTwDr1}

In order to prove Theorem \ref{TwDr1} we need some technical lemmas.

\begin{Lemma} \label{Cor01}
Let $ Y \in \mathbf{S}(\alpha,\kappa,c_Y) $ and $ \hat{Y}(t) = Y(t^{\kappa/\alpha}) $. Then $ \hat{Y} \in \mathbf{S}(\kappa,\kappa,c_Y (\kappa/\alpha)^\kappa) $ and there exist finite and positive
\begin{eqnarray*}
c_1 & = & \inf_{x \in [0,1)} \frac{V_Y(x^{\kappa/\alpha})}{|1-x|^\kappa} \; = \; \inf_{x \in [0,1)} \frac{V_{\hat{Y}}(x)}{|1-x|^\kappa} \; ; \\[1ex]
c_2 & = & \sup_{x \in [0,1)} \frac{V_Y(x^{\kappa/\alpha})}{|1-x|^\kappa} \; = \; \sup_{x \in [0,1)} \frac{V_{\hat{Y}}(x)}{|1-x|^\kappa} \; .
\end{eqnarray*}
Moreover for all $ t,s \geq 0 $,
\[ c_1 |t-s|^{\kappa} \; \leq \; V_{Y}(t^{\kappa/\alpha},s^{\kappa/\alpha}) \; = \; V_{\hat{Y}}(t,s) \; \leq \; c_2 |t-s|^{\kappa} \; . \] \vspace{2ex}
\end{Lemma}
 \textbf{Proof:} \\
Observe that
$ \hat{Y} \in \mathbf{S}(\kappa,\kappa,c_Y (\kappa/\alpha)^\kappa) $
with $ V_{\hat{Y}}(t,s) = V_{Y}(t^{\kappa/\alpha},s^{\kappa/\alpha}) $
and $ v_{\hat{Y}}(x) = v_{Y}(x^{\kappa/\alpha}) $.
Consider function $ f(x) = \frac{v_{\hat{Y}}(x)}{|1-x|^\kappa} $ for $ x \in [0,1) $.
Due to \textbf{S2}, $ \lim_{x \to 1^-} f(x) = c_{\hat{Y}} > 0 $, $ f(0) = 1 $
and
$ f(x) = 0 $ only for $ x = 1 $.
Hence, $ c_1, c_2 > 0 $ introduced in the thesis of Lemma \ref{Cor01} exist.

Moreover, for all $ t \geq s \geq 0 $,
\[ c_1 |t-s|^{\kappa} \; = \; c_1 t^\kappa |1-s/t|^{\kappa} \; \leq \; t^\kappa v_{\hat{Y}}(s/t) \; \leq \; c_2 t^\kappa |1-s/t|^{\kappa} \; = \; c_2 |t-s|^{\kappa} \; . \]
This completes the proof. \hfill \raisebox{-1.5ex}{$ \square $} \\[1ex]

\begin{Lemma} \label{lemTh01}
Let $ \hat{Y} \in \mathbf{S}(\kappa,\kappa,c_{\hat{Y}}) $. For any $ \epsilon \to 0^+ $ there exists $ \delta_\epsilon \to 0^+ $, such that for any $ T > 0 $ and $ A \geq T/\delta_\epsilon $,
\[ (1-\epsilon) c_{\hat{Y}} |t-s|^\kappa \; \leq \; V_{\hat{Y}}(A+t,A+s) \; \leq \; (1+\epsilon) c_{\hat{Y}} |t-s|^\kappa \; , \]
for all $ t,s \in [0,T] $. \\
\end{Lemma}

\noindent \textbf{Proof:} \\
Let $ \epsilon \in (0,1) $ be sufficiently small such that
\begin{equation} \label{lab02} (1-\epsilon) c_{\hat{Y}} |h|^\kappa \; \leq \; V_{\hat{Y}}(1,1-h) \; \leq \; (1+\epsilon) c_{\hat{Y}} |h|^\kappa \; , \end{equation}
for all $ h \in [0,\delta_\epsilon] $ and $ \delta_\epsilon \in (0,1) $ (due to \textbf{S2}).
Then, for any $ T > 0 $, $ A \geq T/\delta_\epsilon $ and $ 0 \leq s \leq t \leq T $ we have
$ \frac{t-s}{A+t} \leq \frac{T}{A} \leq \delta_\epsilon $. Combining the fact  that
\[ V_{\hat{Y}}(A+s,A+t) \; = \; (A+t)^\kappa V_{\hat{Y}}  \left( 1,\frac{A+s}{A+t} \right) \; = \; (A+t)^\kappa V_{\hat{Y}}
\left(1, 1 - \frac{t-s}{A+t} \right) \]
with (\ref{lab02}), for $ h = \frac{t-s}{A+t} \leq \delta_\epsilon $, we obtain the thesis. This completes the proof. \hfill \raisebox{-1.5ex}{$ \square $}

\begin{Lemma} \label{Lem26}
Let $ \hat{Y} \in \mathbf{S}(\kappa,\kappa,c_{\hat{Y}}) $ and $ c_1, c_2 $
be such that the thesis of Lemma \ref{Cor01} holds. Consider a centered Gaussian process $ \{ \overline{X}(t) : t \geq 0 \} $ with $ R_{\overline{X}}(t,s) = \exp(-aV_{\hat{Y}}(t,s)) $ for $ a > 0 $.
\begin{itemize}
\item[\emph{(i)}] Let $ \{ \overline{X}_1(t) : t \geq 0 \} $ and $ \{ \overline{X}_2(t) : t \geq 0 \} $ be centered stationary Gaussian processes with $ R_{\overline{X}_i}(t,s) = \exp(-ac_i V_{B_\kappa}(t,s)) $ respectively. Then for all $ u > 0 $ and any $ T, A=A(u) > 0 $,
\begin{eqnarray*}
\mathbf{P}  \left( \sup_{t \in [0,T]u^{-2/\kappa}} \overline{X}_1(t) > u \right) & \leq & \mathbf{P}  \left( \sup_{t \in [A,A+T]u^{-2/\kappa}} \overline{X}(t) > u \right) \\[1ex]
& \leq & \mathbf{P}  \left( \sup_{t \in [0,T]u^{-2/\kappa}} \overline{X}_2(t) > u \right) .
\end{eqnarray*}
\item[\emph{(ii)}] For any $ \epsilon \to 0^+ $ there exists $ \delta_\epsilon \to 0^+ $, such that for any $ T > 0 $ and $ A=A(u) \geq T/\delta $,
\begin{eqnarray*}
\mathbf{P}  \left( \sup_{t \in [0,T]u^{-2/\kappa}} \overline{X}_1(t) > u \right) & \leq & \mathbf{P}  \left( \sup_{t \in [A,A+T]u^{-2/\kappa}} \overline{X}(t) > u \right) \\[1ex]
& \leq & \mathbf{P}  \left( \sup_{t \in [0,T]u^{-2/\kappa}} \overline{X}_2(t) > u \right) ,
\end{eqnarray*}
where $ \{ \overline{X}_1(t) : t \geq 0 \} $ and $ \{ \overline{X}_2(t) : t \geq 0 \} $ are centered stationary Gaussian processes with $ R_{\overline{X}_i}(t,s) = \exp \left(-a (1+(-1)^i \epsilon) c_{\hat{Y}} V_{B_\kappa}(t,s)\right) $ respectively. \\
\end{itemize}
\end{Lemma}

\noindent \textbf{Proof:} \\
Ad (i): The proof is based on the same argument as given in the proof of Theorem \ref{TwDr2}.
From Lemma \ref{Cor01} we have that for all $ t,s \geq 0 $,
\[ V_{X_1}(t,s) \; \leq \; V_{X}(t,s) \; \leq \; V_{X_2}(t,s) \]
and hence, due to Slepian's inequality, for all $ u > 0 $,
\begin{eqnarray*}
\mathbf{P}  \left( \sup_{t \in [A,A+T]u^{-2/\kappa}} \overline{X}_2(t) > u \right) & \leq & \mathbf{P}  \left( \sup_{t \in [A,A+T]u^{-2/\kappa}} \overline{X}(t) > u \right) \\[1ex]
& \leq & \mathbf{P}  \left( \sup_{t \in [A,A+T]u^{-2/\kappa}} \overline{X}_2(t) > u \right) .
\end{eqnarray*}
Due to stationarity of $ \overline{X}_i(\cdot) $ we obtain the thesis.

Ad (ii): From Lemma \ref{lemTh01}, for any $ \epsilon \to 0^+ $ there exists $ \delta_\epsilon \to 0^+ $, such that for any $ T > 0 $ and $ A \geq T/\delta_\epsilon $,
\[ (1-\epsilon) c_{\hat{Y}} |t-s|^\kappa \; \leq \; V_{\hat{Y}}(t,s) \; \leq \; (1+\epsilon) c_{\hat{Y}} |t-s|^\kappa \; , \]
for all $ t,s \in [A,A+T] $. The same argument as given in the proof of part (i)
completes the proof. \hfill \raisebox{-1.5ex}{$ \square $} \\[1ex]

\begin{Lemma} \label{Lem25}
Suppose that $ \lim_{u \to \infty} f(u)/u = c $, for some $ c > 0 $. Under the notation of Lemma \ref{Lem26}, there exist absolute constants $ F, G > 0 $ such that
\[ \mathbf{P}  \left( \sup_{t \in [A,A+T]u^{-2/\kappa}} \overline{X}(t) > f(u) , \sup_{t \in [t_0,t_0+T]u^{-2/\kappa}} \overline{X}(t) > f(u) \right) \]
\[ \leq \; F T^2 \exp \left( - G (t_0-(A+T))^\kappa \right) \Psi(f(u)) \; , \]
for all $ t_0 > A + T > 0 $, $T\ge 1$ and any $ u \geq u_0 = (2ac_2)^2 (t_0 + T)^{\kappa/2} $.
%, i.e. $ (t_0 + T)u^{-2/\kappa} \leq (2ac_2)^{-1/\kappa} $.
\end{Lemma}

\noindent \textbf{Proof:} \\
The proof follows by argument similar to the one given in ,e.g., Lemma 6.2 in \cite{Deb02} or Theorem 2.1 in \cite{DHL17}.
Thus we present only main steps of the proof.\\
Let $ u_0 = (2ac_2)^2 (t_0 + T)^{\kappa/2} $ and $ \left\{ Z_u(t_1,t_2) : (t_1,t_2) \in [A,A+T]\times[t_0,t_0+T] \right\} $, where $ Z_u(t_1,t_2) = \overline{X}(t_1u^{-2/\kappa}) + \overline{X}(t_2u^{-2/\kappa}) $. Note that
\[ \mathbf{P}  \left( \sup_{t \in [A,A+T]u^{-2/\kappa}} \overline{X}(t) > f(u) , \sup_{t \in [t_0,t_0+T]u^{-2/\kappa}} \overline{X}(t) > f(u) \right) \]
\begin{equation} \label{lab34} \leq \; \mathbf{P}  \left( \sup_{(t_1,t_2) \in [A,A+T]\times[t_0,t_0+T]} Z_u(t_1,t_2) > 2f(u) \right) . \end{equation}
Since $ (t_0 + T)u^{-2/\kappa} \leq (2ac_2)^{-1/\kappa} $ then (from  Lemma \ref{Cor01}) for all $ t_1,t_2 \leq t_0+T $,
\begin{eqnarray}
ac_1 u^{-2}|t_2-t_1|^{\kappa} & \leq & aV_{\hat{Y}} \! \left(t_1u^{-2/\kappa},t_2u^{-2/\kappa}\right) \\[1ex]
\label{lab49} & \leq & ac_2 u^{-2}|t_2-t_1|^{\kappa} \; \leq \; ac_2 \left|(t_0+T)u^{-2/\kappa}\right|^{\kappa} \; \leq \; 1/2 \; .
\end{eqnarray}
Hence, using the fact that $ x \leq 2(1-e^{-x}) \leq  \left(1-e^{-4x}\right) $ for $ x \in [0,1/2] $, we obtain
\begin{eqnarray}
\nonumber V_{\overline{X}} \! \left(t_1u^{-2/\kappa},t_2u^{-2/\kappa}\right) & = & 2  \left( 1 - \exp \left(-aV_{\hat{Y}} \! \left(t_1u^{-2/\kappa},t_2u^{-2/\kappa}\right)\right) \right) \\[1ex]
\label{lab33} & \geq & aV_{\hat{Y}} \! \left(t_1u^{-2/\kappa},t_2u^{-2/\kappa}\right) \; \geq \; ac_1 u^{-2}|t_2-t_1|^{\kappa} \; ; \\[1ex]
\nonumber V_{\overline{X}} \! \left(t_1u^{-2/\kappa},t_2u^{-2/\kappa}\right) & \leq & 2  \left( 1 - \exp \left(-ac_2 u^{-2}|t_2-t_1|^{\kappa}\right) \right) \\[1ex]
\label{lab37} & \leq &  \left( 1 - \exp \left(-4ac_2 u^{-2}|t_2-t_1|^{\kappa}\right) \right)
\end{eqnarray}
for all $ t_1,t_2 \leq t_0+T $. Since
\begin{eqnarray*}
\sigma_{Z_u}^2(t_1,t_2) & = & 2 + 2 \exp \left(-aV_{\hat{Y}} \! \left(t_1u^{-2/\kappa},t_2u^{-2/\kappa}\right)\right) \\[1ex]
& = & 4 - 2  \left( 1 - \exp \left(-aV_{\hat{Y}} \! \left(t_1u^{-2/\kappa},t_2u^{-2/\kappa}\right)\right) \right) ,
\end{eqnarray*}
then from (\ref{lab33}), for any $ (t_1,t_2) \in [A,A+T]\times[t_0,t_0+T] $,
\begin{equation} \label{lab35} 2 \; \leq \; \sigma_{Z_u}^2(t_1,t_2) \; \leq \; 4 - ac_1 u^{-2}(t_0-(A+T))^{\kappa} \; . \end{equation}

Now observe that
\begin{eqnarray*}
\lefteqn{
\mathbf{P}  \left( \sup_{(t,s) \in [A,A+T]\times[t_0,t_0+T]} Z_u(t_1,t_2) > 2f(u) \right) \le}\\
\label{lab39} &\leq & \mathbf{P}  \left( \sup_{(t_1,t_2) \in [A,A+T]\times[t_0,t_0+T]}
\overline{Z}_u(t_1,t_2) > \frac{2f(u)}{\sqrt{ 4 - ac_1 u^{-2}(t_0-(A+T))^{\kappa}}} \right) .
\end{eqnarray*}
Note that  for any $ (t_1,t_2), (s_1,s_2) \in [A,A+T]\times[t_0,t_0+T] $, we have
\begin{eqnarray}
\mathbf{Var}(\overline{Z}_u(t_1,t_2)-\overline{Z}_u(s_1,s_2))
&\leq&
\frac{\mathbf{Var}({Z}_u(t_1,t_2)-{Z}_u(s_1,s_2))}{\sigma_{Z_u}(t_1,t_2) \sigma_{Z_u}(s_1,s_2)}\nonumber\\
&\leq&\frac{1}{2} \: \mathbf{E}  \left(  \left( X(t_1u^{-2/\kappa}) - X(s_1u^{-2/\kappa}) \right) +
\left( X(t_2u^{-2/\kappa}) - X(s_2u^{-2/\kappa}) \right) \right)^2\nonumber\\
 \label{lab36} &\leq& V_X \! \left(t_1u^{-2/\kappa},s_1u^{-2/\kappa}\right) +
 V_X \! \left(t_2u^{-2/\kappa},s_2u^{-2/\kappa}\right)\\
\label{lab38} &\leq& \left( 1 - \exp \left(-4ac_2 u^{-2}|t_1-s_1|^{\kappa}\right) \right) +
\left( 1 - \exp \left(-4ac_2 u^{-2}|t_2-s_2|^{\kappa}\right) \right) ,
\end{eqnarray}
where (\ref{lab36}) follows from inequality $ (x+y)^2 \leq 2(x^2+y^2) $ and (\ref{lab38}) follows from (\ref{lab37}).

Denote $ u^* = \frac{2f(u)}{\sqrt{ 4 - ac_1 u^{-2}(t_0-(A+T))^{\kappa}}} $ and let $ \underline{c}, \overline{c} > 0 $
be constants such that $ \underline{c} \leq \frac{f(u)}{u} \leq \overline{c} $ for all $ u \geq u_0 $. Note (by (\ref{lab49}))
that $ f(u) \leq u^* \leq \sqrt{8/7} f(u) $ for $ u \geq u_0 $.
Hence, $ \underline{c} u \leq u^* \leq \sqrt{8/7} \overline{c} u $ for $ u \geq u_0 $,
and therefore $ u^{-2} \leq \frac{8}{7} \overline{c}^2 (u^*)^{-2} $ for $ u \geq u_0 $.

Consider two independent, identically distributed centered stationary Gaussian processes\\
$ \left\{ Z_{1,u^*}(t_1) : t_1 \geq 0 \right\} $, $ \left\{ Z_{2,u^*}(t_2) : t_2
\geq 0 \right\} $ with $ R_{Z_{1,u^*}}(t_1,s_1) =
\exp \left(-\frac{32}{7}a\overline{c}^2c_2 (u^*)^{-2}|t_1-s_1|^{\kappa}\right) $
and let $ \underline{Z}_{u^*}(t_1,t_2) = \frac{1}{\sqrt{2}}  \left( Z_{1,u^*}(t_1) + Z_{2,u^*}(t_2) \right) $.

Hence, by (\ref{lab38}), for any $ (t_1,t_2), (s_1,s_2) \in [A,A+T]\times[t_0,t_0+T] $,
\begin{eqnarray*}
\lefteqn{ \mathbf{Var}(\overline{Z}_u(t_1,t_2)-\overline{Z}_u(s_1,s_2)) \le}\\
&\leq&  \left( 1 - \exp \left(-4ac_2 u^{-2}|t_1-s_1|^{\kappa}\right) \right) + \left( 1 - \exp \left(-4ac_2 u^{-2}|t_2-s_2|^{\kappa}\right) \right)\\
& \leq & \left( 1 - \exp \left(-\frac{32}{7}a\overline{c}^2c_2 (u^*)^{-2}|t_1-s_1|^{\kappa}\right) \right)  + \left( 1 - \exp \left(-\frac{32}{7}a\overline{c}^2c_2 (u^*)^{-2}|t_2-s_2|^{\kappa}\right) \right)\\
& =&
\mathbf{Var}(\underline{Z}_{u^*}(t_1,t_2)-\underline{Z}_{u^*}(s_1,s_2))
\end{eqnarray*}
and due to Slepian's inequality, we obtain that
\begin{eqnarray}
\mathbf{P}  \left( \sup_{(t_1,t_2) \in [A,A+T]\times[t_0,t_0+T]} \overline{Z}_u(t_1,t_2) > u^* \right)
&\leq&
\mathbf{P}  \left( \sup_{(t_1,t_2) \in [A,A+T]\times[t_0,t_0+T]} \underline{Z}_{u^*}(t_1,t_2) > u^* \right)\nonumber\\
\label{lab40}
&=&  \mathbf{P}  \left( \sup_{(t_1,t_2) \in [0,T]^2} \underline{Z}_{u^*}(t_1,t_2) > u^* \right)%\\
%&=&
%\left( \mathcal{H}_{B_\kappa} \! \left((2ac^2c_2)^{1/\kappa}T\right) \right)^2 \Psi(u^*) (1+o(1)) , \label{lab141}
\end{eqnarray}
as $u\to\infty$,
where equality (\ref{lab40}) follows from stationarity of $ \underline{Z}_u^*(\cdot,\cdot) $.
Now
\begin{eqnarray}
\label{lab50} \lim_{u^*\to\infty} \frac{\mathbf{P} \left( \sup_{(t_1,t_2) \in [0,T]^2} \underline{Z}_{u^*}(t_1,t_2) > u^* \right)}{\Psi(u^*)} & = & \left( \mathcal{H}_{B_\kappa} \left((16a{\overline{c}}^2c_2/7)^{1/\kappa}T\right) \right)^2 \\
\label{lab51} & \leq & \left(\mathcal{H}_{B_\kappa} (1)\right)^2 \max\left(1,\left(16a{\overline{c}}^2c_2/7\right)^{2/\kappa} \right)T^2 ,
\end{eqnarray}
where equality (\ref{lab50}) follows from, e.g., Theorem 2.1 in \cite{Deb02} (see also Theorem 3.1 in \cite{DHL17})
and inequality (\ref{lab51}) follows from the fact that $ \mathcal{H}_{B_\kappa}(AT) \leq T\max(1,A) \mathcal{H}_{B_\kappa}(1) $
for any $ T > 1 $ and $A>0$ (Corollary D.1 in \cite{Pit96}).
Hence, there exists a constant $F'$ (that does not depend on $t_0,A,T$), such that
\begin{eqnarray*}
\mathbf{P} \left( \sup_{(t_1,t_2) \in [0,T]^2} \underline{Z}_{u^*}(t_1,t_2) > u^* \right) \leq F' T^2 \Psi(u^*)
\end{eqnarray*}
\noindent holds for all $ u^* \geq u_0^* = \underline{c}u_0 $ (i.e. $ u \geq u_0 $). Combination of the above with (\ref{lab34}) and (\ref{lab39}) gives
\begin{equation} \label{lab41} \mathbf{P} \left( \sup_{t \in [A,A+T]u^{-2/\kappa}} \overline{X}(t) > f(u) , \sup_{t \in [t_0,t_0+T]u^{-2/\kappa}} \overline{X}(t) > f(u) \right) \leq F' T^2 \Psi(u^*) , \end{equation}
\noindent for $ u \geq u_0 $.

Since (using inequality $\frac{1}{1-x} \geq 1+x $, for $ x \geq 0 $)
\begin{eqnarray*}
(u^*)^2 & = & \frac{4f^2(u)}{4 - ac_1 u^{-2}(t_0-(A+T))^{\kappa}} \; \geq \; f^2(u) + \frac{ac_1}{4}  \left( \frac{f(u)}{u} \right)^2 (t_0-(A+T))^{\kappa} \\[1ex]
& \geq & f^2(u) + \frac{ac_1\underline{c}^2}{4} (t_0-(A+T))^{\kappa} \; ,
\end{eqnarray*}
then,
\begin{eqnarray}
\nonumber \Psi(u^*) & \leq & \frac{\exp  \left( - \frac{1}{2}  \left( f^2(u) + \frac{ac_1\underline{c}^2}{4} (t_0-(A+T))^{\kappa} \right) \right)}{\sqrt{2 \pi} \sqrt{f^2(u) + \frac{ac_1\underline{c}^2}{4} (t_0-(A+T))^{\kappa}}} \\[1ex]
\label{lab42} & \leq & \frac{\exp  \left( - \frac{1}{2} f^2(u) \right)}{\sqrt{2 \pi} f(u)} \exp  \left( - \frac{ac_1\underline{c}^2}{8} (t_0-(A+T))^{\kappa} \right) .
\end{eqnarray}
Combination of (\ref{lab41}) with (\ref{lab42}) gives
\begin{eqnarray*}
\lefteqn{\mathbf{P}  \left( \sup_{t \in [A,A+T]u^{-2/\kappa}} \overline{X}(t) > f(u) ,
\sup_{t \in [t_0,t_0+T]u^{-2/\kappa}} \overline{X}(t) > f(u) \right) \le}\\
&\leq&
 F_1 F' T^2 \exp  \left( - \frac{ac_1\underline{c}^2}{8} (t_0-(A+T))^{\kappa} \right) \Psi(f(u)) \; ,
 \end{eqnarray*}
for any $ u \geq u_0 $ and some positive constant $ F_1 $, such that $ \frac{\exp \left( - \frac{1}{2} f^2(u) \right)}{\sqrt{2 \pi} f(u)} \leq F_1 \Psi(f(u)) $ for $ u > u_0 $.
This completes the proof with $ F = F_1 F' $ and $ G = \frac{ac_1\underline{c}^2}{8} $. \hfill \raisebox{-1.5ex}{$ \square $}

\begin{Lemma} \label{Lem27}
Under the notation of Lemma \ref{Lem26}, there exist absolute constants $ F, G > 0 $ such that
\[ \mathbf{P}  \left( \sup_{t \in [A,A+T]u^{-2/\kappa}} \overline{X}(t) > u , \sup_{t \in [A+T,A+2T]u^{-2/\kappa}} \overline{X}(t) > u \right) \]
\[ \leq \; F  \left( T^2 \exp \left( - G \sqrt{T^\kappa} \right) + \sqrt{T} \right) \Psi(u) \; , \]
for all $ A>0, T > 1 $ and any $ u \geq u_0 = (2ac_2)^2 (A+2T)^{\kappa/2} $, i.e. $ (A+2T)u^{-2/\kappa} \leq (2ac_2)^{-1/\kappa} $.
\\
\end{Lemma}

\noindent \textbf{Proof:} \\
Let $ u_0 = (2ac_2)^2 (A+2T)^{\kappa/2} $ and $ \overline{X}_u(t) = \overline{X}(tu^{-2/\kappa}) $. We have
\begin{eqnarray*}
\lefteqn{\mathbf{P}  \left( \sup_{t \in [A,A+T]} \overline{X}_u(t) > u , \sup_{t \in [A+T,A+2T]} \overline{X}_u(t) > u \right)}\\
&=& \mathbf{P}  \left( \sup_{t \in [A,A+T]} \overline{X}_u(t) > u , \left\{ \sup_{t \in [A+T,A+T\sqrt{T}]} \overline{X}_u(t) > u \vee \sup_{t \in [A+T+\sqrt{T},A+2T]} \overline{X}_u(t) > u \right\} \right)\\
&\leq&
\mathbf{P}  \left( \sup_{t \in [A,A+T]} \overline{X}_u(t) > u , \sup_{t \in [A+T+\sqrt{T},A+2T+\sqrt{T}]} \overline{X}_u(t) > u \right)\\
&& +
\mathbf{P}  \left( \sup_{t \in [A+T,A+T+\sqrt{T}]} \overline{X}_u(t) > u \right) \\
&\leq&
F_1 T^2 \exp \left( - G \sqrt{T^\kappa} \right) \Psi(u) \; + \; \mathbf{P}  \left( \sup_{t \in [A+T,A+T+\sqrt{T}]} \overline{X}_u(t) > u \right) ,
\end{eqnarray*}
where the last inequality follows from Lemma \ref{Lem25} with $ t_0 = A + T + \sqrt{T} $.
Applying Lemma \ref{Lem26} (i) and Proposition \ref{Cor09}  to the above,
we obtain that for sufficiently large $ u \geq u_0 $,
\begin{eqnarray}
\lefteqn{\nonumber
\mathbf{P}  \left( \sup_{t \in [A,A+T]} \overline{X}_u(t) > u , \sup_{t \in [A+T,A+2T]} \overline{X}_u(t) > u \right)}\\
& \leq & \nonumber
F_1 T^2 \exp \left( - G \sqrt{T^\kappa} \right) \Psi(u) + \mathcal{H}_{B_\kappa} \! \left(a^{1/\kappa}\sqrt{T}\right) \Psi(u) (1+o(1))\\
\label{lab151} & \leq &
F_1 T^2 \exp \left( - G \sqrt{T^\kappa} \right) \Psi(u) + \max\left(1,a^{1/\kappa}\right) \mathcal{H}_{B_\kappa}(1) \sqrt{T} \Psi(u) (1+o(1))\\
&\leq &
F  \left( T^2 \exp \left( - G \sqrt{T^\kappa} \right) + \sqrt{T} \right) \Psi(u)\nonumber
\end{eqnarray}
for some constant $ F > 0 $, where (\ref{lab151}) follows from subadditivity of $ \mathcal{H}_{B_\kappa}(\cdot) $
(Corollary D.1 in \cite{Pit96}).
This completes the proof. \hfill \raisebox{-1.5ex}{$ \square $}
\\

\noindent \textbf{Proof of Theorem \ref{TwDr1}:} \\
First, we prove the thesis for $ \hat{Y} \in \mathbf{S}(\kappa,\kappa,c_Y (\kappa/\alpha)^\kappa) $, where $ \hat{Y}(t) = Y(t^{\kappa/\alpha}) $.
Let $ c_1, c_2 $ be constants such that the thesis of Lemma \ref{Cor01} holds.
Consider a centered Gaussian process $ \{ \overline{X}(t) : t \geq 0 \} $
with $ R_{\overline{X}}(t,s) = \exp(-V_{\hat{Y}}(t,s)) $ and let
$ \overline{X}_u(t) = \overline{X}(tu^{-2/\kappa}) $. Let $ n \in \mathbb{N}$ and choose
$ \epsilon_n \in (0,1) $, $ \delta_{\epsilon_n} = 1/n $ in such a way that the thesis of Lemma \ref{Lem26} (ii) holds.
We find a lower and an upper bound separately. \\

\noindent \textbf{Upper bound:} \\
Let $ T \in \mathbb{N} $ be such that $ T > n $. Due to Bonferroni's inequality, for any $ u > 0 $,
\[ \mathbf{P}  \left( \sup_{t \in \left[0,T^2\right]} \overline{X}_u(t) > u \right) \; \leq \; \mathbf{P}  \left( \sup_{t \in [0,nT]} \overline{X}_u(t) > u \right) \; + \; \sum_{k=n}^{T-1} \mathbf{P}  \left( \sup_{t \in [kT,(k+1)T]} \overline{X}_u(t) > u \right) . \]
Applying Lemma \ref{Lem26} to the right side of the inequality above, by Proposition \ref{Cor09}, we obtain that
\[ \mathcal{H}_{\hat{Y}}(T^2) \; \leq \; \mathcal{H}_{B_\kappa} \! \left(c_2^{1/\kappa}nT\right) \; + \; (T-n) \mathcal{H}_{B_\kappa} \! \left( \left((1+\epsilon_n)c_{\hat{Y}}\right)^{1/\kappa}T\right) \]
and hence,
\begin{equation} \label{lab114} \frac{\mathcal{H}_{\hat{Y}}(T^2)}{T^2} \; \leq \; \frac{\mathcal{H}_{B_\kappa} \! \left(c_2^{1/\kappa}nT\right)}{T^2} \; + \; \frac{ \left((1+\epsilon_n)c_{\hat{Y}}\right)^{1/\kappa}T^2}{T^2} \frac{\mathcal{H}_{B_\kappa} \! \left( \left((1+\epsilon_n)c_{\hat{Y}}\right)^{1/\kappa}T\right)}{ \left((1+\epsilon_n)c_{\hat{Y}}\right)^{1/\kappa}T} \; . \end{equation}
Since $ \lim_{S \to \infty} S^{-1} \mathcal{H}_{B_\kappa} (S) = \mathcal{H}_{B_\kappa} $ then after sending $ T \to \infty $ in (\ref{lab114}), we get
\[ \limsup_{T \to \infty} \frac{\mathcal{H}_{\hat{Y}}(T)}{T} \; \leq \;  \left((1+\epsilon_n)c_{\hat{Y}}\right)^{1/\kappa} \mathcal{H}_{B_\kappa} \; . \]
Since the bound obtained above holds for any $ \epsilon_n \to 0^+ $, then
\begin{equation} \label{lab43} \limsup_{T \to \infty} \frac{\mathcal{H}_{\hat{Y}}(T)}{T} \; \leq \;  \left(c_{\hat{Y}}\right)^{1/\kappa} \mathcal{H}_{B_\kappa}. \end{equation}

\noindent \textbf{Lower bound:} \\
For $ T \in \mathbb{N} $ such that $ T > n $, with $ \Delta_{k} = [kT,(k+1)T] $, again from Bonferroni's inequality observe that for any $ u > 0 $,
\[ \mathbf{P}  \left( \sup_{t \in \left[0,T^2\right]} \overline{X}_u(t) > u \right) \; \geq \; \mathbf{P}  \left( \sup_{t \in \left[nT,T^2\right]} \overline{X}_u(t) > u \right) \]
\[ \geq \; \sum_{k=n}^{T-1} \mathbf{P}  \left( \sup_{t \in \Delta_{k}} \overline{X}_u(t) > u \right) \; - \; \sum_{1 \leq k \leq l}^{T-1} \mathbf{P}  \left( \sup_{t \in \Delta_{k}} \overline{X}_u(t) > u , \sup_{t \in \Delta_{l}} \overline{X}_u(t) > u \right) \]
\begin{equation} \label{lab44} \geq \; \sum_{k=n}^{T-1} \mathbf{P}  \left( \sup_{t \in \Delta_{k}} \overline{X}_u(t) > u \right) \; - \; \Sigma_1 - \Sigma_2 \; , \end{equation}
where
\begin{eqnarray*}
\Sigma_1 & = & \sum_{k=1}^{T-2} \mathbf{P}  \left( \sup_{t \in \Delta_{k}} \overline{X}_u(t) > u , \sup_{t \in \Delta_{k+1}} \overline{X}_u(t) > u \right) , \\[1ex]
\Sigma_2 & = & \sum_{1 \leq k < l \neq k+1}^{T-1} \mathbf{P}  \left( \sup_{t \in \Delta_{k}} \overline{X}_u(t) > u , \sup_{t \in \Delta_{l}} \overline{X}_u(t) > u \right) .
\end{eqnarray*}
By Lemma \ref{Lem26} (ii) and Proposition \ref{Cor09}, as $ u \to \infty $, we estimate
\begin{equation} \label{lab45} \sum_{k=n}^{T-1} \mathbf{P}  \left( \sup_{t \in \Delta_{k}} \overline{X}_u(t) > u \right) \; \geq \; (T-n) \mathcal{H}_{B_\kappa} \! \left( \left((1-\epsilon_n)c_{\hat{Y}}\right)^{1/\kappa}T\right) \Psi(u) (1+o(1)) \; . \end{equation}
From Lemma \ref{Lem27}, for sufficiently large $ u $, we have
\begin{equation} \label{lab46} \Sigma_1 \; \leq \; T F_1  \left( T^2 \exp \left( - G_1 \sqrt{T^\kappa} \right) + \sqrt{T} \right) \Psi(u) \; . \end{equation}
From Lemma \ref{Lem25}, for sufficiently large $ u $, we have
\begin{equation} \label{lab47} \Sigma_2 \; \leq \; T^2 F_2 T^2 \exp \left( - G_2 T^\kappa \right) \Psi(u) \; . \end{equation}

Applying (\ref{lab45}), (\ref{lab46}) and (\ref{lab47}) to (\ref{lab44}), (and using Proposition \ref{Cor09}),
we obtain
\begin{eqnarray*}
\frac{\mathcal{H}_{\hat{Y}}(T^2)}{T^2} & \geq & \frac{(T-n) \mathcal{H}_{B_\kappa} \! \left( \left((1-\epsilon_n)c_{\hat{Y}}\right)^{1/\kappa}T\right)}{T^2} \\[1ex]
& & - \; \frac{ F_1  \left( T^3 \exp \left( - G_1 \sqrt{T^\kappa} \right) + T^{3/2} \right) + F_2 T^4 \exp \left( - G_2 T^\kappa \right) }{T^2} \; .
\end{eqnarray*}
Sending $ T \to \infty $ and then $ \epsilon_n \to 0^+ $ in the inequality above to get
\begin{equation} \label{lab48} \liminf_{T \to \infty} \frac{\mathcal{H}_{\hat{Y}}(T)}{T} \; \geq \;  \left(c_{\hat{Y}}\right)^{1/\kappa} \mathcal{H}_{B_\kappa} \; . \end{equation}

From upper bound (\ref{lab43}) and lower bound (\ref{lab48}) we conclude that
\[ \lim_{T \to \infty} \frac{\mathcal{H}_{\hat{Y}}(T)}{T} \; = \;  \left(c_{\hat{Y}}\right)^{1/\kappa} \mathcal{H}_{B_\kappa} \; . \]
For $ Y \in \mathbf{S}(\alpha,\kappa,c_Y) $, $ c_{\hat{Y}} = c_Y (\kappa/\alpha)^\kappa $ and from Lemma \ref{Lem23} (ii) it follows that
\[ \frac{\kappa}{\alpha}  \left(c_Y\right)^{1/\kappa} \mathcal{H}_{B_\kappa} \; = \;  \left(c_{\hat{Y}}\right)^{1/\kappa} \mathcal{H}_{B_\kappa} \; = \; \lim_{T \to \infty} \frac{\mathcal{H}_{\hat{Y}}(T)}{T} \; = \; \lim_{T \to \infty} \frac{\mathcal{H}_{Y}(T^{\kappa/\alpha})}{T} \; = \; \lim_{T \to \infty} \frac{\mathcal{H}_{Y}(T)}{T^{\alpha/\kappa}} \; . \]
This completes the proof. \hfill \raisebox{-1.5ex}{$ \square $} \\[1ex]

\section*{Acknowledgments}
Kamil Tabi{\'s} was partially supported by the NCN Project 2011/03/N/ST1/00163 (2012--2014).

\bibliographystyle{plain}%{ieeetr}
\bibliography{loc}

\begin{thebibliography}{10}

\bibitem{Adl90}
R.~Adler.
\newblock An introduction to continuity, extrema, and related topics for
  general {G}aussian processes.
\newblock IMS, 1990.

\bibitem{Albin}
J.M.P. Albin and H.~Choi.
\newblock A new proof of an old result by {P}ickands.
\newblock {\em Electron. Commun. Probab.}, 15:339--345, 2010.

\bibitem{DHJ18}
L.~Bai, K.~D{\c{e}}bicki, E.~Hashorva, and L.~Luo.
\newblock On generalised {P}iterbarg constants.
\newblock {\em Methodology and Computing in Applied Probability},
  20(1):137--164, 2018.

\bibitem{Bojdecki}
T.~Bojdecki, L.~G Gorostiza, and A.~Talarczyk.
\newblock Sub-fractional {B}rownian motion and its relation to occupation
  times.
\newblock {\em Statistics \& Probability Letters}, 69(4):405--419, 2004.

\bibitem{Mic04}
K.~Burnecki and Z.~Michna.
\newblock Simulation of {P}ickands constants.
\newblock {\em Probability and Mathematical Statistics}, 22:193--199, 2002.

\bibitem{DEH17}
K.~D\c{e}bicki, S.~Engelke, and E.~Hashorva.
\newblock Generalized {P}ickands constants and stationary max-stable processes.
\newblock {\em Extremes}, 20(3):493--517, 2017.

\bibitem{DeH16}
K.~D\c{e}bicki and E.~Hashorva.
\newblock On extremal index of max-stable processes.
\newblock {\em Probability and Mathematical Statistics}, 27(2):299--317, 2017.

\bibitem{Deb02}
K~D{\c{e}}bicki.
\newblock Ruin probability for {G}aussian integrated processes.
\newblock {\em Stochastic Processes and their Applications}, 98(1):151--174,
  2002.

\bibitem{DeH19}
K.~D{\c{e}}bicki and E.~Hashorva.
\newblock Approximation of supremum of max-stable stationary processes \&
  {P}ickands constants.
\newblock {\em Submitted for publication}, 2019.

\bibitem{DHJR18}
K.~D{\c{e}}bicki, E.~Hashorva, L.~Ji, and T.~Rolski.
\newblock Extremal behavior of hitting a cone by correlated {B}rownian motion
  with drift.
\newblock {\em Stochastic Processes and their Applications},
  128(12):4171--4206, 2018.

\bibitem{DHT1}
K.~D{\c{e}}bicki, E.~Hashorva, L.~Ji, and K.~Tabi{\'s}.
\newblock Extremes of vector-valued {G}aussian processes: {E}xact asymptotics.
\newblock {\em Stochastic Processes and their Applications},
  125(11):4039--4065, 2015.

\bibitem{DHL17}
K.~D{\c{e}}bicki, E.~Hashorva, and P.~Liu.
\newblock Uniform tail approximation of homogenous functionals of {G}aussian
  fields.
\newblock {\em Advances in Applied Probability}, 49(4):1037--1066, 2017.

\bibitem{DebMan}
K.~D{\c{e}}bicki and M.~Mandjes.
\newblock Exact overflow asymptotics for queues with many {G}aussian inputs.
\newblock {\em Journal of Applied Probability}, 40(3):704--720, 2003.

\bibitem{DMR03}
K.~D{\c{e}}bicki, Z.~Michna, and T.~Rolski.
\newblock Simulation of the asymptotic constant in some fluid models.
\newblock {\em Stochastic Models}, 19(3):407--423, 2003.

\bibitem{DebTab}
K.~D{\c{e}}bicki and K.~Tabi{\'s}.
\newblock Extremes of the time-average of stationary {G}aussian processes.
\newblock {\em Stochastic Processes and their Applications}, 121(9):2049--2063,
  2011.

\bibitem{DiM15}
A.~B. Dieker and T.~Mikosch.
\newblock Exact simulation of {B}rown-{R}esnick random fields at a finite
  number of locations.
\newblock {\em Extremes}, 18:301--314, 2015.

\bibitem{DzZ04}
K.~Dzhaparidze and H.~Van~Zanten.
\newblock A series expansion of fractional {B}rownian motion.
\newblock {\em Probability theory and related fields}, 130(1):39--55, 2004.

\bibitem{HKP19}
E~Hashorva, S~Kobelkov, and V.I. Piterbarg.
\newblock On maximum of {G}aussian process with unique maximum point of its
  variance.
\newblock {\em arXiv:1901.09753}, 2019.

\bibitem{Ho}
Ch. Houdr{\'e} and J.~Villa.
\newblock An example of infinite dimensional quasi-helix.
\newblock {\em Contemporary Mathematics}, 336:195--202, 2003.

\bibitem{KonsPit}
D.G. Konstant and V.I. Piterbarg.
\newblock Extreme values of the cyclostationary {G}aussian random process.
\newblock {\em Journal of Applied Probability}, 30(1):82--97, 1993.

\bibitem{Lam}
J.~Lamperti.
\newblock Semi-stable stochastic processes.
\newblock {\em Transactions of the American mathematical Society},
  104(1):62--78, 1962.

\bibitem{Lei}
P.~Lei and D.~Nualart.
\newblock A decomposition of the bifractional {B}rownian motion and some
  applications.
\newblock {\em Statistics \& Probability Letters}, 79(5):619--624, 2009.

\bibitem{LiShao}
W.~Li and Q.~Shao.
\newblock Lower tail probabilities for {G}aussian processes.
\newblock {\em The Annals of Probability}, 32(1A):216--242, 2004.

\bibitem{Mic17}
Z~Michna.
\newblock Remarks on {P}ickands' theorem.
\newblock {\em Probability and Mathematical Statistics}, (37):373--393, 2017.

\bibitem{Pic69}
J.~Pickands.
\newblock Asymptotic properties of the maximum in a stationary {G}aussian
  process.
\newblock {\em Transactions of the American Mathematical Society}, 145:75--86,
  1969.

\bibitem{Pit96}
V.~I. Piterbarg.
\newblock {\em Asymptotic methods in the theory of {G}aussian processes and
  fields}, volume 148.
\newblock American Mathematical Soc., 2012.

\bibitem{PiP79}
V.I. Piterbarg and V.P. Prisiazhniuk.
\newblock Asymptotic analysis of the probability of large excursions for a
  nonstationary {G}aussian process.
\newblock {\em Teoriia Veroiatnostei i Matematicheskaia Statistika},
  (18):121--134, 1978.

\bibitem{Wit96}
R.~A. Vitale.
\newblock The {W}ills functional and {G}aussian processes.
\newblock {\em Ann. Probab.}, 24(4):2172--2178, 1996.

\end{thebibliography}

\end{document}